\documentclass[12pt]{article}

\usepackage{amsmath,amsthm,amsfonts,amssymb,amscd}
\usepackage{tikz}
\usepackage{color}
\usepackage{cite, url}

\newtheorem{theorem}{Theorem}
\newtheorem{conjecture}[theorem]{Conjecture}

\newtheorem{problem}[theorem]{Problem}

\begin{document}

\title{On the local metric dimension of $K_4$-free graphs}

\author{
Ali Ghalavand $^{a,}$\thanks{Corresponding author, email: \texttt{alighalavand@nankai.edu.cn}}
\and
Sandi Klav\v zar $^{b,c,d}$
\and 
Xueliang Li $^{a}$
}

\maketitle

\begin{center}
$^a$ Center for Combinatorics, Nankai University, Tianjin 300071, China \\
\medskip

$^b$ Faculty of Mathematics and Physics,  University of Ljubljana, Slovenia\\
\medskip

$^b$ Institute of Mathematics, Physics and Mechanics, Ljubljana, Slovenia\\
\medskip

$^d$ Faculty of Natural Sciences and Mathematics,  University of Maribor, Slovenia\\
\end{center}

\begin{abstract}
Let \( G \) be a graph of order \( n(G) \), local metric dimension \( \dim_l(G) \), and clique number \( \omega(G) \). It has been conjectured that if \( n(G) \geq \omega(G) + 1 \geq 4 \), then \( \dim_l(G) \leq \left( \frac{\omega(G) - 2}{\omega(G) - 1} \right) n(G) \). In this paper the conjecture is confirmed for the case \( \omega(G) = 3 \). Consequently, a problem regarding the local metric dimension of planar graphs is also resolved.
\end{abstract}


\medskip\noindent
\textbf{AMS Math.\ Subj.\ Class.\ (2020)}: 05C12, 05C69

\section{Introduction}

Let \( G \) be a simple connected graph with vertex set \( V(G) \) and edge set \( E(G) \). We denote the clique number of \( G \) by \( \omega(G) \) and the order of \( G \) by \( n(G) \). The open neighborhood \( N_G(u) \) of \( u\in V(G) \) is the set of vertices adjacent to \( u \). The degree \( d_G(u) \) of \( u \) is the cardinality of \( N_G(u) \). The maximum degree of \( G \) is denoted by \( \Delta(G) \). 
If \( V' \) is a subset of the vertex set \( V(G) \), then the notation \( G[V'] \) refers to the subgraph of \( G \) that is induced by \( V' \). Additionally, if \( H \) and \( H' \) are subgraphs of \( G \), then \( E_G(H, H') \) denotes the set of edges connecting vertices from \( H \) to vertices in \( H' \).
The union \( G \cup G' \) of graphs \( G \) and \( G' \) is the graph with vertex set \( V(G) \cup V(G') \) and edge set \( E(G) \cup E(G') \).
For a positive integer \( l \geq 1 \), the notation \( [l] \) represents the set \( \{ 1, \ldots, l \} \). We also set \( [0] = \emptyset \).

The distance \( d_G(u, v) \) between vertices \( u, v\in V(G) \) is the length of shortest $u,v$-path in \( G \). The vertices $u$ and $v$ are {\em distinguished} by $w\in V(G)$, or equivalently, $w$ {\em distinguishes} $u$ and $v$, if $d_G(u,w)\neq d_G(v,w)$. A subset \( W \subseteq V(G) \) is a {\em resolving set} for \( G \) if, for any two vertices \( u \) and \( v \) in \( V(G) - W \), there exists at least one vertex in \( W \) that distinguishes \( u \) and \( v \). Similarly, \( W \) is referred to as a {\em local resolving set} if, for any adjacent vertices \( u \) and \( v \) from \( V(G) - W \), there is a vertex in \( W \) that distinguishes \( u \) and \( v \). The {\em metric dimension} \( \dim(G) \) and the {\em local metric dimension} \( \dim_l(G) \) of \( G \) are defined as the sizes of the smallest resolving sets and the smallest local resolving sets for \( G \), respectively. Clearly,  \( \dim_l(G) \leq \dim(G) \). 

The metric dimension of graphs has a rich history, initially defined by Harary and Melter~\cite{13} and by Slater~\cite{25}. Determining the metric dimension is known to be NP-complete for general graphs~\cite{19} and also for restricted cases involving planar graphs with a maximum degree six~\cite{6}. Research in this area is extensive, partly because metric dimension has numerous real-world applications, including robot navigation, image processing, privacy in social networks, and tracking intruders in networks. The 2023 overview~\cite{survey2} of the essential results and applications of metric dimension contains well over 200 references.

Several variations of metric dimension gained a wider attention. The survey~\cite{survey1} focusing on these variants also cites over 200 papers. One particularly interesting variation is the local metric dimension, introduced in 2010 by Okamoto, Phinezy, and Zhang~\cite{Okamoto1}. Like the standard metric dimension, the local metric dimension is computationally challenging \cite{9,10} and has been explored in several studies \cite{Abrishami1, 3, 4, fitriani, Ghalavand1, 17, lal-2023, rodriguez-2016}. We should also mention closely related research on the fractional local metric dimension \cite{javaid-2024} and the nonlocal metric dimension~\cite{klavzar-2023}. Okamoto et al.~\cite{Okamoto1} proved several important relationships between local metric dimension and clique number: $\dim_l(G) = n(G) - 1$ if and only if $G \cong K_{n(G)}$; $\dim_l(G) = n(G) - 2$ if and only if $\omega(G) = n(G) - 1$; $\dim_l(G) = 1$ if and only if $G$ is bipartite; and $\dim_l(G) \geq \max\left \{ \lceil \log_2 \omega(G)\rceil, n(G) - 2^{n(G) - \omega(G)}\right\}$. Furthermore, Abrishami et al.~\cite{Abrishami1} established that $\dim_l(G) \leq \frac{2}{5}n(G)$ when $\omega(G) = 2$ and $n(G) \geq 3$. They also posed:

\begin{problem} {\rm  \cite[Problem 1]{Abrishami1}} 
\label{con1}
If $G$ is a planar graph with $n(G) \geq 2$, is it then true that  
\[\dim_l(G) \leq \left\lceil\frac{n(G)+1}{2}\right\rceil\,?\]
\end{problem}

In~\cite{Abrishami1} it was confirmed that Problem~\ref{con1} has a positive answer for triangle-free planar graphs. On the other hand, it was demonstrated in~\cite{Ghalavand2} that the problem has a negative answer when \( \omega(G) = 4 \). Additionally, in~\cite{Ghalavand1} it was proved that \( \dim_l(G) \leq \left( \frac{\omega(G) - 1}{\omega(G)} \right) n(G) \), with equality holding if and only if \( G \cong K_{n(G)} \), thus verifying a conjecture from~\cite{Abrishami1}. 

Our second main motivation is:

\begin{conjecture} {\rm \cite[Conjecture 2]{Ghalavand1}}
\label{con2}
If $G$ is a graph with $n(G) \geq \omega(G)+1 \geq 4$, then 
\[\dim_l(G)\leq\left(\frac{\omega(G)-2}{\omega(G)-1}\right)n(G)\,.\]
\end{conjecture}

It is demonstrated in~\cite{Ghalavand1} that if Conjecture~\ref{con2} is true, then the bound is asymptotically best possible, and that Conjecture~\ref{con2} holds for all graphs $G$ with \(\omega(G) \in \{n(G)-1, n(G)-2, n(G)-3\}\). It is also established that when \(\omega(G) = n(G)-2\), then \(n(G)-4 \leq \dim_l(G) \leq n(G)-3\), and that when \(\omega(G) = n(G)-3\), then \(n(G)-8 \leq \dim_l(G) \leq n(G)-3\).

The main result of this paper reads as follows. 

\begin{theorem}\label{th3}
If $G$ is a graph with $n(G) \geq 4>\omega(G)$, then 
\[\dim_l(G) \leq \left\lfloor \frac{n(G)}{2} \right\rfloor.\]
\end{theorem}

Theorem~\ref{th3} implies that Problem~\ref{con1} and Conjecture~\ref{con2} have positive answers when \(\omega(G) = 3\). Additionally, there exist infinitely many planar graphs that achieve the equality stated in the theorem. For any positive odd number \(n\), let \(G = \frac{n-1}{2}K_2 + K_1\), that is, a graph constructed from \(\frac{n-1}{2}\) disjoint complete graphs \(K_2\) by adding a new vertex and connecting it to all the vertices of the \(\frac{n-1}{2}K_2\). It is straightforward to observe that the local metric dimension of \(G\) equals \(\lfloor \frac{n}{2} \rfloor\). 
Therefore, there are infinitely many planar graphs \(G\) such that \(\dim_l(G) = \left\lfloor \frac{n(G)}{2} \right\rfloor\). 

So, Problem~\ref{con1} has a positive answer for planar graph $G$ with $\omega(G)\le 3$, while we know from before that there are planar graph $G$ containing $K_4$ for which the bound of Problem~\ref{con1} does not apply. It is important to highlight that Conjecture~\ref{con2} remains unresolved for graphs \( G \) where \( 4 \leq \omega(G) \leq n(G) - 4 \).

\section{Proof of Theorem~\ref{th3}}
\label{sec:proof}

\begin{figure}[ht!]
\begin{center}
\begin{tikzpicture}
\clip(-6.5,0) rectangle (4,4.5);
\draw [line width=1.2pt] (-5,4)-- (-6.02,3.53);
\draw [line width=1.2pt] (-6.02,3.53)-- (-5,3);
\draw [line width=1.2pt] (-5,4)-- (-5,3);
\draw [line width=1.2pt] (-5,4)-- (-4,3.48);
\draw [line width=1.2pt] (-5,3)-- (-4,3.48);
\draw [line width=1.2pt] (-3,4)-- (-3,3);
\draw [line width=1.2pt] (-3,4)-- (-2,4);
\draw [line width=1.2pt] (-2,4)-- (-3,3);
\draw [line width=1.2pt] (-2,4)-- (-2,3);
\draw [line width=1.2pt] (-1,4)-- (-1,3);
\draw [line width=1.2pt] (-1,3)-- (0,3);
\draw [line width=1.2pt] (0,3)-- (-1,4);
\draw [line width=1.2pt] (1,4)-- (1,3);
\draw [line width=1.2pt] (1,3)-- (2,3);
\draw [line width=1.2pt] (2,3)-- (2,4);
\draw [line width=1.2pt] (2,4)-- (1,4);
\draw [line width=1.2pt] (-6,2)-- (-6,1);
\draw [line width=1.2pt] (-6,2)-- (-5,2);
\draw [line width=1.2pt] (-5,2)-- (-4.99,0.96);
\draw [line width=1.2pt] (-3,2)-- (-4,2);
\draw [line width=1.2pt] (-4,2)-- (-4,1);
\draw [line width=1.2pt] (-4,2)-- (-3,1);
\draw [line width=1.2pt] (-2,2)-- (-1,2);
\draw [line width=1.2pt] (-2,1)-- (-1,1);
\draw [line width=1.2pt] (0,2)-- (0,1);
\draw [line width=1.2pt] (0,2)-- (1,2);
\draw [line width=1.2pt] (2,1)-- (3,1);
\draw (-5.27,2.9) node[anchor=north west] {$F_1$};
\draw (-2.72,2.9) node[anchor=north west] {$F_2$};
\draw (-0.77,2.9) node[anchor=north west] {$F_3$};
\draw (1.28,2.9) node[anchor=north west] {$F_4$};
\draw (-5.76,0.8) node[anchor=north west] {$F_5$};
\draw (-3.75,0.8) node[anchor=north west] {$F_6$};
\draw (-1.76,0.8) node[anchor=north west] {$F_7$};
\draw (0.28,0.8) node[anchor=north west] {$F_8$};
\draw (2.26,0.8) node[anchor=north west] {$F_9$};
\begin{scriptsize}
\fill [color=black] (-5,4) circle (2.0pt);
\draw[color=black] (-4.92,4.21) node {$a_1$};
\fill [color=black] (-6.02,3.53) circle (2.0pt);
\draw[color=black] (-6,3.74) node {$a_4$};
\fill [color=black] (-5,3) circle (2.0pt);
\draw[color=black] (-5.2,3.31) node {$a_3$};
\fill [color=black] (-4,3.48) circle (2.0pt);
\draw[color=black] (-3.9,3.7) node {$a_2$};
\fill [color=black] (-3,4) circle (2.0pt);
\draw[color=black] (-2.91,4.21) node {$b_1$};
\fill [color=black] (-3,3) circle (2.0pt);
\draw[color=black] (-3.2,3.21) node {$b_4$};
\fill [color=black] (-2,4) circle (2.0pt);
\draw[color=black] (-1.9,4.21) node {$b_2$};
\fill [color=black] (-2,3) circle (2.0pt);
\draw[color=black] (-1.8,3.21) node {$b_3$};
\fill [color=black] (-1,4) circle (2.0pt);
\draw[color=black] (-0.94,4.21) node {$c_1$};
\fill [color=black] (-1,3) circle (2.0pt);
\draw[color=black] (-1.25,3.1) node {$c_3$};
\fill [color=black] (0,3) circle (2.0pt);
\draw[color=black] (0.22,3.1) node {$c_2$};
\fill [color=black] (1,4) circle (2.0pt);
\fill [color=black] (1,3) circle (2.0pt);
\fill [color=black] (2,3) circle (2.0pt);
\fill [color=black] (2,4) circle (2.0pt);
\fill [color=black] (-6,2) circle (2.0pt);
\fill [color=black] (-6,1) circle (2.0pt);
\fill [color=black] (-5,2) circle (2.0pt);
\fill [color=black] (-4.99,0.96) circle (2.0pt);
\fill [color=black] (-3,2) circle (2.0pt);
\fill [color=black] (-4,2) circle (2.0pt);
\fill [color=black] (-4,1) circle (2.0pt);
\fill [color=black] (-3,1) circle (2.0pt);
\fill [color=black] (-2,2) circle (2.0pt);
\fill [color=black] (-1,2) circle (2.0pt);
\fill [color=black] (-2,1) circle (2.0pt);
\fill [color=black] (-1,1) circle (2.0pt);
\fill [color=black] (0,2) circle (2.0pt);
\fill [color=black] (0,1) circle (2.0pt);
\fill [color=black] (1,2) circle (2.0pt);
\fill [color=black] (2,1) circle (2.0pt);
\fill [color=black] (3,1) circle (2.0pt);
\end{scriptsize}
\end{tikzpicture}
\caption{ All non-complete graphs with at most $4$ vertices and no isolated vertex.}
\label{fig1}
\end{center}
\end{figure}
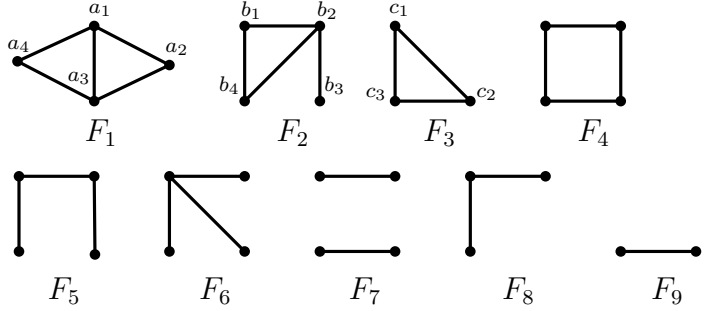

We begin the proof by describing a key approach to it. First, let \( F_i \), \( i \in [9] \), be the graphs of order at most $4$ and with no isolated vertices, except $K_4$, as illustrated in Fig.~\ref{fig1}. 

Let now $G$ be a graph with $n(G) \geq 4>\omega(G)$. In the following, we will sequentially select in $G$ and in its induced subgraphs maximum sets of vertex disjoint induced subgraphs isomorphic to some $F_i$. Such a selection is not necessarily unique, but we will select one and fix it, so the following notation is well-defined for the purposes of the proof.     

\begin{itemize}
\item Let $\mathcal{F}_1(G)$ be a maximum set of vertex disjoint induced subgraphs of $G$ isomorphic to $F_1$. 
\item Set $G_1 = G$. For $i=2,3,\ldots,9$, let $\mathcal{F}_i(G)$ be a maximum set of vertex disjoint induced subgraphs of $G_{i}=G_{i-1} - \bigcup_{H \in \mathcal{F}_{i-1}(G)}V(H)$ isomorphic to $F_i$.
\item Note that $G_{9}-\bigcup_{H \in \mathcal{F}_{9}(G)}V(H)$ is a set of isolated vertices, let $\mathcal{F}_{10}(G)$ be the set of these induced subgraphs isomorphic to $K_1$.
\item Let \( \mathcal{F}(G) = \{\mathcal{F}_{i}(G):\ i \in [10]\} \), and call it a {\em local vertex division} of \( G \).
\item For \( i \in [10] \), let  \(V_i= \cup_{H \in \mathcal{F}_{i}(G)} V(H) \).
\end{itemize}

The sets \(V_i \), \( i \in [10] \), form a partition of \( V(G) \). We further emphasize that some of these sets may be empty, that is, some of the sets \( \mathcal{F}_i(G) \) may be empty. The definition of the local vertex division and the assumption $\omega(G)\leq 3$ imply the following facts, where we use vertex labels from Fig.~\ref{fig1}. 

\begin{enumerate}
  \item Let $H\in\mathcal{F}_1(G)$ and $v\in V(G-H)$. If for some $i\in\{2,4\}$, $va_i\in E(G)$, then there exists $j\in\{1,3\}$ such that $va_j\not\in E(G)$.
  \item  Let $H\in\mathcal{F}_2(G)$ and $v\in V(G_2-H)$. Then at most one vertex in the set $\{b_1,b_2,b_4\}$ is adjacent to $v$.
  \item  Let $H\in\mathcal{F}_4(G)$ and $v\in V(G_4-H)$. Then at most two vertices of $H$ are  adjacent to $v$, and $G[V(H)\cup\{v\}]$ does not contain cycles of length 3.
  \item  Let $H\in\mathcal{F}_5(G)$ and $v\in V(G_5-H)$. Then at most two vertices of $H$ can be adjacent to $v$, and $G[V(H) \cup \{v\}]$ does not contain cycles of length 3 or 4.
  \item For $i\in\{6,7,8,9\}$, let $H\in\mathcal{F}_i(G)$ and $v\in V(G_i-H)$. Then at most one vertex  of $H$ is  adjacent to $v$.
\end{enumerate}

Based on the above five statements we claim that there exists a subset $S$ of $V(G-\cup_{H \in \mathcal{F}_3(G)}H)$ such that 
\begin{itemize}
\item for any $H \in \mathcal{F}_i(G)$, $i\in [10]\setminus \{3\}$, we have $|S\cap V(H)|\geq\frac{1}{2}n(H)$, and 
\item $V(G)-S$ is a local resolving set for $G$. 
\end{itemize}
We construct $S$ as follows, where we initially set $S = \emptyset$. For any graph $H_1\in \mathcal{F}_1(G)$, we add to $S$ its vertices $a_2$ and $a_4$ (see Fig.~\ref{fig1}). Similarly, for any $H_2\in \mathcal{F}_2(G)$, add to $S$ either $b_1$ and $b_2$, or $b_2$ and $b_4$. Next, for any graph $H\in \mathcal{F}_4(G)\cup\mathcal{F}_5(G)\cup\mathcal{F}_6(G)\cup\mathcal{F}_7(G)$, add to $S$ two arbitrary non-adjacent vertices of $H$.  Further, for any $H_8\in \mathcal{F}_8(G)$, add to $S$ two arbitrary vertices of $H_8$, while for any element $H_9\in \mathcal{F}_9(G)$, add to $S$ an arbitrary vertex of $H_9$. Finally, set $S=S\cup\mathcal{F}_{10}(G)$. Since it is evident that $S$ satisfies the required two conditions, the claim is proved. 

If $\mathcal{F}_3(G)$ is empty, then the above constructed local resolving set 
$V(G)-S$ yields the theorem's assertion. Hence assume in the rest that $\mathcal{F}_3(G) \neq \emptyset$. Then the following statements hold, where in each case the argument is based on the fact that otherwise we could increase the cardinality of $\mathcal{F}_1(G)$ or $\mathcal{F}_2(G)$. 

\begin{enumerate}
\item[I.] Let \( H \) and \( H' \) be two disjoint elements in the set \( \mathcal{F}_3(G) \). Then, for any \( h \in V(H) \) and any \( h' \in V(H') \), it holds that \( hh' \not\in E(G) \).
  \item[II.] There is no edge between the vertices that lie on the elements of $\mathcal{F}_3(G)$ and those that lie on the elements of $\cup_{i=4}^{10}\mathcal{F}_i(G)$. In other words, for each $u \in \left(\cup_{H \in \mathcal{F}_3(G)} V(H)\right)$ and $v \in \left(\cup_{i=4}^{10} \cup_{H \in \mathcal{F}_i(G)} V(H)\right)$, it holds that $uv \notin E(G)$.
 \item[III.] Let $F\in\mathcal{F}_3(G)$ and $V(F)=\{c_1,c_2,c_3\}$. If for $v\in\left(\cup_{i=1}^2\cup_{H \in \mathcal{F}_i(G)} V(H)\right)$ and $i\in[3]$, we have $vc_{i}\in\,E(G)$, then $v$ distinguishes  either $(c_i$ and $c_j)$ or $(c_i$ and $c_k)$, where $[3]-\{i\}=\{j,k\}$.
 \item[IV.] If for $H \in \mathcal{F}_2(G)$ and $F,F'\in\mathcal{F}_3(G)$, we have $V(H)\cap\left(\cup_{v\in\,V(F)}N_G(v) \right)\neq\emptyset$ and $V(H)\cap\left(\cup_{u\in\,V(F')}N_G(u) \right)\neq\emptyset$, then 
 \[|\left(V(H)\cap\left(\cup_{v\in\,V(F)}N_G(v) \right)\right)\cup\left(V(H)\cap\left(\cup_{u\in\,V(F')}N_G(u) \right)\right)|=1.\]
\end{enumerate}

In the rest of the proof we are going to add to the above constructed set $S$ some of the vertices from the triangles from $\mathcal{F}_3(G)$, such that $V(G) - S$ remains a local resolving set for $G$. For this sake, some additional notation is needed. 

\begin{itemize}
\item Let $\mathcal{F}_3(G)\neq\emptyset$, and let $\mathcal{A}$ be a subset of $\mathcal{F}_3(G)$.  For $H \in \mathcal{F}_1(G)$, let $\eta_1(H, \mathcal{F}_3(G), \mathcal{A})$ represent the set of elements in $\mathcal{F}_3(G) - \mathcal{A}$ such that at least one vertex from each of these elements is adjacent to some vertices of $H$, that is, 
\[
\eta_1(H, \mathcal{F}_3(G), \mathcal{A}) = \{ F : F \in \mathcal{F}_3(G) - \mathcal{A}, |E_G(H, F)| \geq 1 \}.
\]
\item Analogously, for $H' \in \mathcal{F}_2(G)$, let $\eta_2(H', \mathcal{F}_3(G), \mathcal{A})$ represent the set of elements in $\mathcal{F}_3(G) - \mathcal{A}$ such that at least one vertex from each of these elements is adjacent to some vertices of $H'$, that is, 
\[
\eta_2(H', \mathcal{F}_3(G), \mathcal{A}) = \{ F : F \in \mathcal{F}_3(G) - \mathcal{A}, |E_G(H', F)| \geq 1 \}.
\]
\item Also, let $ H'' \in \mathcal{F}_1(G) \cup \mathcal{F}_2(G)$. For a subset $ U $ of $V(H'')$ and a subset $Y $ of $\mathcal{F}_3(G)$, the notation $D(U, Y)$ represents the largest set of two-subsets $\{x, y\}$ such that the following conditions hold: (i) $ xy $ is an edge from a triangle in $ Y $; (ii) $ x $ and $y $ are distinguished by a vertex from $ U $; (iii) any triangle from $Y$ has at most  one two-subset in \( D(U, Y) \).
\end{itemize}

\medskip\noindent
\underline{$1^{\rm st}$ process:}
\begin{enumerate}
\item[(1.1)] Set $\mathcal{A}=\emptyset$ and consider the above-defined set $S$.
\item[(1.2)] Select $H\in \mathcal{F}_1(G)$ such that $\eta_1(H, \mathcal{F}_3(G), \mathcal{A})$ has maximum cardinality.  
\item[(1.3)] If  $|\eta_1(H, \mathcal{F}_3(G), \mathcal{A})|\leq1$, then return  $S$ and $\mathcal{A}$,  and end the process,  otherwise go to (1.4).
\item[(1.4)] If  $|\eta_1(H, \mathcal{F}_3(G), \mathcal{A})|\geq4$, then set 
\begin{align*}
S & = \left(S-V(H)\right)\cup\,D(V(H), \eta_1(H, \mathcal{F}_3(G), \mathcal{A})),\\
\mathcal{A} & = \mathcal{A}\cup\eta_1(H, \mathcal{F}_3(G), \mathcal{A}),
\end{align*}
and proceed to (1.2), otherwise go to (1.5).
\item[(1.5)] If  $|\eta_1(H, \mathcal{F}_3(G), \mathcal{A})| \in \{2,3\}$ and $h\in V(H)$ has the property  $\eta_1(H, \mathcal{F}_3(G), \mathcal{A})=\eta_1(H-h, \mathcal{F}_3(G), \mathcal{A})$, then set
\begin{align*}
S & = \left(S-V(H)\right)\cup\,D(V(H)-\{h\}, \eta_1(H, \mathcal{F}_3(G), \mathcal{A}))\cup\{h\}, \\
\mathcal{A} & = \mathcal{A}\cup\eta_1(H, \mathcal{F}_3(G), \mathcal{A}),
\end{align*} 
and proceed to (1.2).
\end{enumerate}

\newpage\noindent
\underline{$2^{\rm nd}$ process:} 
\begin{enumerate} 
\item[(2.1)] Consider the sets $\mathcal{A}$ and $S$ that are returned in the $1^{\rm st}$ process.
\item[(2.2)] Select $H\in \mathcal{F}_2(G)$ such that $\eta_2(H, \mathcal{F}_3(G), \mathcal{A})$ has maximum cardinality.  
\item[(2.3)] If  $|\eta_2(H, \mathcal{F}_3(G), \mathcal{A})|\leq1$, then return  $S$ and $\mathcal{A}$,  and end the process, otherwise go to (2.4).
\item[(2.4)] For $h\in V(H)$ such that $d_H(h)\geq 2$ and $h$ is adjacent to no vertex from $\eta_2(H, \mathcal{F}_3(G), \mathcal{A})$, set
\begin{align*}
S & = \left(S-V(H)\right)\cup\{h\}\cup\,D(V(H)-\{h\}, \eta_2(H, \mathcal{F}_3(G), \mathcal{A})), \\
\mathcal{A} & = \mathcal{A}\cup\eta_2(H, \mathcal{F}_3(G), \mathcal{A}),
\end{align*}
and proceed to (2.2).
\end{enumerate}

\medskip\noindent
\underline{$3^{\rm rd}$ process:} 
\begin{enumerate} 
\item[(3.1)] Consider the sets $\mathcal{A}$ and $S$ that are returned in the $2^{\rm nd}$ process. 
\item[(3.2)] If $\mathcal{F}_3(G)-\mathcal{A}\neq\emptyset$, then go to (3.3), otherwise  return  $S$ and $\mathcal{A}$,  and end the process.
\item[(3.3)] Take an element $F$ of $ \mathcal{F}_3(G)$.  If there exists an element $H$ in $\mathcal{F}_1(G)$ such that for a vertex $h$ of it, $d_H(h)=3$ and $h$ is adjacent to a vertex of $F$, then set 
\begin{align*}
S & = S\cup D(\{h\},\{F\}),\\
\mathcal{A} & = \mathcal{A}\cup\{F\}, 
\end{align*}
and proceed to (3.2), otherwise go to (3.4).
\item[(3.4)] If there exists an element $H$ in $\mathcal{F}_2(G)$ such that for a vertex $h$ of it, $d_{H}(h)\leq2$ and $h$ is adjacent to a vertex of $F$, then for $h_1,h_2\in V(H)$ such that $d_{H}(h_1)=2$, $d_{H}(h_2)=3$, and $h_1\neq h$, set
\begin{align*}
S & = (S-V(H))\cup\{h_1,h_2\}\cup\,D(\{h\},\{F\}),\\
\mathcal{A} & = \mathcal{A}\cup\{F\}, 
\end{align*}
and proceed to (3.2).
\end{enumerate}

Before starting the $4^{\rm th}$ process, let's consider the sets $\mathcal{A}$ and $S$ that are returned in the $3^{\rm rd}$ process. Also, if $z_1$ and $z_2$ are integers such that $z_1+z_2=|\mathcal{F}_3(G)-\mathcal{A}|$, then let $\mathcal{F}_3(G)-\mathcal{A}=\{F_i^1:i\in[z_1]\}\cup\{F_j^2:j\in[z_2]\}$, where for $i\in[z_1]$, $F_i^1$ has adjacency with an element in $\mathcal{F}_1(G)$, and for $j\in[z_2]$, $F_j^2$ has adjacency with an element in $\mathcal{F}_2(G)$. Now,  assume for $i\in[z_1]$, $H_i^1\in\mathcal{F}_1(G)$ has adjacency with $F_i^1$, and for $j\in[z_2]$,  $H_j^2\in\mathcal{F}_2(G)$ has adjacency with $F_j^2$. Also, for $i\in[z_1]$, let's consider $V(F_i^1)=\{f^1_{i_1},f^1_{i_2},f^1_{i_3}\}$, $V(H_i^1)=\{h^1_{i_1},h^1_{i_2},h^1_{i_3},h^1_{i_4}\}$, $d_{H_i^1}(h^1_{i_2})=d_{H_i^1}(h^1_{i_4})=2$, $d_{H_i^1}(h^1_{i_1})=d_{H_i^1}(h^1_{i_3})=3$,  $f^1_{i_1}h^1_{i_2}\not\in\,E(G)$, and $f^1_{i_2}h^1_{i_2}\in\,E(G)$. Plus, for $j\in[z_2]$, let's consider $V(F_j^2)=\{f^2_{j_1},f^2_{j_2},f^2_{j_3}\}$, $V(H_j^2)=\{h^2_{j_1},h^2_{j_2},h^2_{j_3},h^2_{j_4}\}$, $d_{H_j^2}(h^2_{j_1})=d_{H_j^2}(h^2_{j_4})=2$, $d_{H_j^2}(h^2_{j_2})= 3$, $d_{H_j^2}(h^2_{j_3})= 1$, and $f^2_{j_2}h^2_{j_2}\in\,E(G)$.

\medskip\noindent
\underline{$4^{\rm th}$ process:}
\begin{enumerate}
\item[(4.1)] Consider the set $S$ that is returned in the $3^{\rm rd}$ process and set $Z_1=[z_1]$.
\item[(4.2)] If there is $i\in\,Z_1$ and $x\in\mathcal{F}_{10}(G)$ such that $\{xh^1_{i_1},xh^1_{i_4}\}\subseteq\,E(G)$,  then go to (4.3), otherwise  return  $S$,  $Z_1$, $\mathcal{F}_{10}(G)$  and end the process.
\item[(4.3)] Set 
\begin{align*}
S & = (S-\{x,h^1_{i_2}\})\cup\{h^1_{i_3},f^1_{i_1},f^1_{i_2}\}, \\
Z_1 & = Z_1-\{i\}, \\
\mathcal{F}_{10}(G) & = \mathcal{F}_{10}(G)-\{x\}, 
\end{align*}
and proceed to (4.2).
\end{enumerate}

Before starting the $5^{\rm th}$ process, let's set up $\mathcal{F}_{9}(G)=\{F^9_1,\ldots,F^9_{|\mathcal{F}_{9}(G)|}\}$.

\medskip\noindent
\underline{$5^{\rm th}$ process:}
\begin{enumerate}
\item[(5.1)] Consider the sets $S$ and $Z_1$  that are returned in the $4^{\rm th}$ process. Also, for $i\in[|\mathcal{F}_{9}(G)|]$, set  $V^9_i=V(F^9_{i})$.
\item[(5.2)] If there is $i\in\,Z_1$ and $j\in[|\mathcal{F}_{9}(G)|]$,  such that for a vertex $x$ in  $V^9_j$ we have $\{xh^1_{i_1},xh^1_{i_4}\}\subseteq\,E(G)$,  then go to (5.3), otherwise  return  $S$,  $Z_1$,   and end the process.
\item[(5.3)] Set 
\begin{align*}
S & = ((S\cup V^9_j)-\{x,h^1_{i_2}\})\cup\{h^1_{i_3},f^1_{i_1},f^1_{i_2}\}, \\  Z_1 & = Z_1-\{i\}, \\ 
V^9_j & = V^9_j-\{x\}, 
\end{align*}
and proceed to (5.2).
\end{enumerate}

Before starting the $6^{\rm th}$ process, let's set up $\mathcal{F}_{8}(G)=\{F^8_1,\ldots,F^8_{|\mathcal{F}_{8}(G)|}\}$.

\medskip\noindent
\underline{$6^{\rm th}$ process:}
\begin{enumerate}
\item[(6.1)] Consider the sets $S$ and $Z_1$  that are returned in the $5^{\rm th}$ process. Also, for $i\in[|\mathcal{F}_{8}(G)|]$, set  $V^8_i=V(F^8_{i})$.
\item[(6.2)] If there is $i\in\,Z_1$ and $j\in[|\mathcal{F}_{8}(G)|]$,  such that for a vertex $x$ in  $V^8_j$ we have $\{xh^1_{i_1},xh^1_{i_4}\}\subseteq\,E(G)$,  then go to (6.3), otherwise  return  $S$,  $Z_1$,   and end the process.
\item[(6.3)] Set 
\begin{align*}
S & = ((S\cup\,V^8_j)-\{x,h^1_{i_2}\})\cup\{h^1_{i_3},f^1_{i_1},f^1_{i_2}\}, \\  Z_1 & = Z_1-\{i\}, \\ 
V^8_j & = V^8_j-\{x\}, 
\end{align*}
and proceed to (6.2).
\end{enumerate}

We set up next $\mathcal{F}_{7}(G)=\{\{x_{i_1}x_{i_2},y_{i_1}y_{i_2}\}: i\in[|\mathcal{F}_{8}(G)|]\}$.

\medskip\noindent
\underline{$7^{\rm th}$ process:}
\begin{enumerate}
\item[(7.1)] Consider the sets $S$ and $Z_1$  that are returned in the $6^{\rm th}$ process. Also, for $i\in[|\mathcal{F}_{7}(G)|]$ and $a\in\{x,y\}$, set  $V^7_{a_i}=\{a_{i_1},a_{i_2}\}$.
\item[(7.2)] If there is $i\in\,Z_1$ and $j\in[|\mathcal{F}_{7}(G)|]$,  such that for $a\in\{x,y\}$ and $v\in\,V^7_{a_j}$  we have $\{vh^1_{i_1},vh^1_{i_4}\}\subseteq\,E(G)$,  then go to (7.3), otherwise  return  $S$,  $Z_1$,   and end the process.
\item[(7.3)] Set 
\begin{align*}
S & = ((S\cup\,V^7_{a_j})-\{v,h^1_{i_2}\})\cup\{h^1_{i_3},f^1_{i_1},f^1_{i_2}\},\\ 
Z_1 & = Z_1-\{i\}, \\
V^7_{a_j} & = V^7_{a_j}-\{v\}, 
\end{align*}
and proceed to (7.2).
\end{enumerate}

\newpage\noindent
Let's set up $\mathcal{F}_{6}(G)=\{F^6_i: i\in[|\mathcal{F}_{6}(G)|]\}$. Also, for $i\in[|\mathcal{F}_{6}(G)|]$,  set up  $V(F^6_i)=\{f^6_i:i\in[4]\}$ and $d_{F^6_i}(f^6_4)=3$.

\medskip\noindent
\underline{$8^{\rm th}$ process:}
\begin{enumerate}
\item[(8.1)] Consider the sets $S$ and $Z_1$  that are returned in the $7^{\rm th}$  process. Also, for $i\in[|\mathcal{F}_{6}(G)|]$, set  $V^6_i=V(F^6_i)$.
\item[(8.2)] If there is $i\in\,Z_1$ and $j\in[|\mathcal{F}_{6}(G)|]$,  such that for a vertex $v$ in  $V^6_{j}$  we have $\{vh^1_{i_1},vh^1_{i_4}\}\subseteq\,E(G)$,  then go to (8.3), otherwise  return  $S$,  $Z_1$, and end the process.
\item[(8.3)] If $|V^6_{j}|\leq2$, then  set 
\begin{align*}
S & = ((S\cup\,V^6_{j})-\{v,h^1_{i_2}\})\cup\{h^1_{i_3},f^1_{i_1},f^1_{i_2}\},\\  Z_1 & = Z_1-\{i\},\\
V^6_{j} & = V^6_{j}-\{v\}, 
\end{align*}
and proceed to (8.2), otherwise set 
\begin{align*}
S & = ((S\cup\,(V^6_{j}\cap\{f^6_i:i\in[3]\})-\{v,h^1_{i_2}\})\cup\{h^1_{i_3},f^1_{i_1},f^1_{i_2}\}, \\  
Z_1 & = Z_1-\{i\}, \\
V^6_{j} & = V^6_{j}-\{v\}, 
\end{align*}
and proceed to (8.2).
\end{enumerate}

Set up now $\mathcal{F}_{5}(G)=\{F^5_i: i\in[|\mathcal{F}_{5}(G)|]\}$, and for $X\subseteq V(G)$, let $\overline{X}$ be a maximum subset of $X$ such that $E(G[X])=\emptyset$.

\medskip\noindent
\underline{$9^{\rm th}$ process:}
\begin{enumerate}
\item[(9.1)] Consider the sets $S$ and $Z_1$  that are returned in the $8^{\rm th}$ process. Also, for $i\in[|\mathcal{F}_{5}(G)|]$, set  $V^5_i=V(F^5_i)$.
\item[(9.2)] If there is $i\in\,Z_1$ and $j\in[|\mathcal{F}_{5}(G)|]$,  such that for a vertex $v$ in  $V^5_{j}$  we have $\{vh^1_{i_1},vh^1_{i_4}\}\subseteq\,E(G)$,  then go to (9.3), otherwise  return  $S$,  $Z_1$,   and end the process.
\item[(9.3)] Set 
\begin{align*}
X & = V^5_j-\{v\}, \\
S & = ((S\cup\,\overline{X})-\{v,h^1_{i_2}\})\cup\{h^1_{i_3},f^1_{i_1},f^1_{i_2}\},\\ 
Z_1 & = Z_1-\{i\},\\
V^5_{j} & = V^5_{j}-\{v\}, 
\end{align*}
and proceed to (9.2).
\end{enumerate}

\newpage\noindent
Next, set up $\mathcal{F}_{4}(G)=\{F^4_i: i\in[|\mathcal{F}_{4}(G)|]\}$. 

\medskip\noindent
\underline{$10^{\rm th}$ process:}
\begin{enumerate}
\item[(10.1)] Consider the sets $S$ and $Z_1$  that are returned in the $9^{\rm th}$ process. Also, for $i\in[|\mathcal{F}_{4}(G)|]$, set  $V^4_i=V(F^4_i)$.
\item[(10.2)] If there is $i\in\,Z_1$ and $j\in[|\mathcal{F}_{4}(G)|]$,  such that for a vertex $v$ in  $V^4_{j}$  we have $\{vh^1_{i_1},vh^1_{i_4}\}\subseteq\,E(G)$,  then go to (10.3), otherwise  return  $S$,  $Z_1$,   and end the process.
\item[(10.3)] Set 
\begin{align*}
X & = V^4_j-\{v\},\\
S & = ((S\cup\,\overline{X})-\{v,h^1_{i_2}\})\cup\{h^1_{i_3},f^1_{i_1},f^1_{i_2}\},\\  
Z_1 & = Z_1-\{i\}, \\ 
V^4_{j} & = V^4_{j}-\{v\}, 
\end{align*}
and proceed to (10.2).
\end{enumerate}

\medskip\noindent
\underline{$11^{\rm th}$ process:}
\begin{enumerate}
\item[(11.1)] Consider the sets $S$ and $Z_1$ that are returned in the $10^{\rm th}$  process. Also, set   $Z_2=[z_2]$. 
\item[(11.2)] If there are $i\in\,Z_1$ and $j\in\,Z_2$ such that $\{h^1_{i_1}h^2_{j_3},h^1_{i_4}h^2_{j_3}\}\subseteq\,E(G)$,  then go to (11.3), otherwise  return  $S$,  $Z_1$, $Z_2$,  and end the process.
\item[(11.3)] Set 
\begin{align*}
S & = (S-\{h^1_{i_2}\})\cup\{h^1_{i_3},f^1_{i_1},f^1_{i_2},f^1_{j_1}\},\\
Z_1 & = Z_1-\{i\},\\
Z_2 & = Z_2-\{j\},
\end{align*}
and proceed to (11.2).
\end{enumerate}

\medskip\noindent
\underline{$12^{\rm th}$ process:}
\begin{enumerate}
\item[(12.1)] Consider the sets $S$ and $Z_1$ that are returned in the $11^{\rm th}$ process.  
\item[(12.2)] If there are $i,j\in\,Z_1$ such that $\{h^1_{i_1}h^2_{j_4},h^1_{i_4}h^2_{j_4}\}\subseteq\,E(G)$,  then go to (12.3), otherwise  return  $S$,  $Z_1$,   and end the process.
\item[(12.3)] Set 
\begin{align*}
S & = (S-\{h^1_{i_2},h^2_{j_2},h^2_{j_4}\})\cup\{h^1_{i_3},f^1_{i_1},f^1_{i_2},h^2_{j_3},f^2_{j_1},f^2_{j_2}\},\\  
Z_1 & = Z_1-\{i,j\},
\end{align*}
and proceed to (12.2).
\end{enumerate}

\medskip\noindent
\underline{$13^{\rm th}$ process:}
\begin{enumerate}
\item[(13.1)] Consider the sets $S$ and $Z_1$ that are returned in the $12^{\rm th}$ process.  
\item[(13.2)] If  $Z_1\neq\emptyset$, then take an element $i$ of $Z_1$ and go to (13.3), otherwise  return  $S$   and end the process.
\item[(13.3)] Set
\begin{align*}
S & = (S-\{h^1_{i_2}\})\cup\{h^1_{i_3},f^1_{i_1},f^1_{i_2}\},\\  
Z_1 & = Z_1-\{i\},
\end{align*}
and proceed to (13.2).
\end{enumerate}

\medskip\noindent
\underline{$14^{\rm th}$ process:}
\begin{enumerate}
\item[(14.1)] Consider the set $\mathcal{F}_{10}(G)$ that is returned in the $4^{\rm th}$ process, the set $Z_2$ that is returned in the $11^{\rm th}$ process, and the set $S$ that is returned in the $13^{\rm th}$ process.
\item[(14.2)] If there is $i\in\,Z_2$ and $x\in\mathcal{F}_{10}(G)$ such that $\{xh^2_{i_2},xh^2_{i_3}\}\subseteq\,E(G)$,  then go to (14.3), otherwise  return  $S$,  $Z_2$,  and end the process.
\item[(14.3)] Set 
\begin{align*}
S & = (S-\{x,h^2_{i_2}\})\cup\{h^2_{i_3},f^2_{i_1},f^2_{i_2}\},\\  
Z_2 & = Z_2-\{i\},\\
\mathcal{F}_{10}(G) & =\mathcal{F}_{10}(G)-\{x\},
\end{align*}
and proceed to (14.2).
\end{enumerate}

\medskip\noindent
\underline{$15^{\rm th}$ process:}
\begin{enumerate}
\item[(15.1)] Consider the sets $S$ and $Z_2$ that are returned in the $14^{\rm th}$ process.  
\item[(15.2)] If  $Z_2\neq\emptyset$, then take an element $i$ of $Z_2$ and go to (15.3), otherwise  return  $S$   and end the process.
\item[(15.3)] Set 
\begin{align*}
S & = (S-\{h^2_{i_2}\})\cup\{h^2_{i_3},f^2_{i_1},f^2_{i_2}\},\\
Z_2 & = Z_2-\{i\},
\end{align*}
and proceed to (15.2).
\end{enumerate}

Now, let's examine the set \( S \) that is produced in the $15^{\rm th}$ processes. It is clear that \( |S| \geq \frac{n(G)}{2} \). Furthermore, by utilizing $\omega(G) \leq 3$ and the maximality of \( \mathcal{F}_i(G) \) for \( i \in [9] \), we can observe that \( V(G) - S \) serves as a local resolving set for \( G \). Since \( \dim_l(G) \) is an integer, we have proved Theorem~\ref{th3}.

\section*{Acknowledgments}


The research of Ali Ghalavand and Xueliang Li was supported the by NSFC No.\ 12131013 and 12161141006. Sandi Klav\v{z}ar was supported by the Slovenian Research Agency (ARIS) under the grants P1-0297, N1-0355, and N1-0285.

\end{document}